\documentclass[12pt]{article}
\usepackage{amsfonts,amssymb,amsthm,amsmath,latexsym,bm,enumerate,authblk,color,shapepar,mathrsfs}
\usepackage{amsmath}
\usepackage[colorlinks,linktocpage,linkcolor=blue,pdfpagelabels=false]{hyperref}
\usepackage{bookmark}
\usepackage{graphicx}  %
\usepackage{float}  %
\usepackage{subfigure}
\usepackage[graphicx]{realboxes}
\usepackage{makecell} 
\usepackage{soul} 
\usepackage{color, xcolor} 
\allowdisplaybreaks

\numberwithin{equation}{section}
\newtheorem{theorem}{Theorem}[section]
\newtheorem{lemma}[theorem]{Lemma}

\newtheorem{remark}{Remark}[section]

\newtheorem{problem}{Problem}[section]
\newtheorem{algorithm}{Algorithm}[section]
\newtheorem{example}{Example}[section]

\title{Inverse source problem of  variable-exponent time-fractional diffusion equations}

\author[1]{Li Hu\thanks{2211400016@nbu.edu.cn}}
\author[2]{Zhiyuan Li\thanks{Corresponding author: lizhiyuan@nbu.edu.cn}}
\author[3]{Zhidong Zhang\thanks{zhangzhidong@mail.sysu.edu.cn}}
\affil[1,2]{\normalsize{School of Mathematics and Statistics, Ningbo University, 818 Fenghua Road, Ningbo, 315211, Zhejiang, China.}}
\affil[3]{\normalsize{School of Mathematics (Zhuhai), Sun Yat-sen University, Zhuhai 519082, China; Guangdong Province Key Laboratory of Computational Science, Sun Yat-sen University, Guangzhou 510000, China.}}

\date{}

\begin{document}

\maketitle{}

\begin{abstract} 
This paper investigates an inverse temporal source problem for a variable-exponent
time-fractional diffusion equation with time-dependent diffusion coefficients.
By the perturbation method, the model is equivalently transformed into a constant-order
fractional diffusion equation with a convolution perturbation term. Existence, uniqueness,
and energy estimates of strong and weak solutions to the forward problem are established,
and for a separable source term the uniqueness and Lipschitz stability of the temporal
source function are proved using local integral observations. A numerical reconstruction
algorithm combining Tikhonov regularization with the conjugate gradient method is
developed with a convergence error estimate, and numerical experiments verify its
effectiveness and stability for different variable-exponent functions, noise levels and
source functions. Since fractional advection--diffusion models describe the turbulent
transport of airborne pollutants, the method is also applicable in atmospheric science
for recovering the emission history of a pollution source from sparse measurements.

\medskip 

\noindent{\it Keywords:}  variable-exponent fractional equation; inverse source problem; uniqueness; Lipschitz stability; Tikhonov regularization

\medskip

\end{abstract}

\section{Introduction}

Anomalous diffusion phenomena are widely observed in numerous interdisciplinary fields, 
including physics, chemistry, biology, and other fields, as reviewed in \cite{Deng2020Modeling}. 
Unlike the Gaussian diffusion described by classical Brownian motion, the most typical macroscopic feature of anomalous 
diffusion is that the mean square displacement of particles exhibits a nonlinear power-law relationship with time \cite{Chepizhko2013DiffusionSA}, 
namely \(\langle x^2(t)\rangle \sim t^{\alpha}\). Depending on the range of the exponent \(\alpha\), anomalous diffusion can be 
further classified into subdiffusion (\(0<\alpha<1\)) and superdiffusion (\(1<\alpha<2\)). 
Among these, subdiffusion is particularly prevalent and appears prominently in various contexts, 
such as solute transport in porous media \cite{XU2024142693}, drug release in biological tissues \cite{Ghezal2025AnomalousDT}, 
confinement transport in plasma turbulence \cite{2005PhRv94f5003D}, 
and memory effects in stock markets \cite{Lo1991}, all exhibiting significant subdiffusive characteristics.

In a wide range of practical applications, the geometric structure or physical properties of the diffusion 
medium may vary with time or space, leading to a non-constant fractional order. 
For instance, based on mean square displacement data of calcium ion diffusion in dendrites, 
Sun \cite{Sun2009Variableorder,Sun2010} found that 
the sublinear growth exponent is actually a linear function of time \(t\), 
thereby proposing the variable-exponent subdiffusion model. Similarly, in shale gas reservoirs, 
gas transport in nanoscale organic pores deviates from the classical Darcy's law and exhibits significant 
subdiffusive characteristics; fractional models have been widely employed to describe methane transport 
behavior in shale matrices \cite{chang2024study}, 
and the viscoelastic deformation of artificial fractures \cite{peng2020fractional}.
In the field of shape memory polymers(SMPs), Li \cite{Li2017Variable} pointed out that variable-exponent fractional differential equation models 
 are more suitable than constant-exponent ones for describing the shape memory 
 behavior of SMPs under varying temperature histories.
In the study of magnetorheological elastomers(MREs), 
Poojary and Gangadharan \cite{poojary2018integer} found that variable-exponent fractional constitutive models can more accurately predict 
the viscoelastic response of MREs under the combined influence of frequency and magnetic fields. On the mathematical side, 
Zheng \cite{zheng2024local} introduced a variable-exponent to locally modify the subdiffusion equation, which resolves 
the incompatibility between nonlocal operators and local initial conditions, thereby eliminating 
the initial singularity of the solution while preserving the heavy-tailed behavior away from the initial time. 
Patnaik \cite{Patnaik2020} provided a detailed review of the modeling advantages 
of variable-exponent fractional operators in the aforementioned complex 
systems; further applications can be found in \cite{song2021variable,sun2019review}. 
Consequently, the variable-exponent subdiffusion model demonstrates great potential in both practical applications 
and theoretical research.

Atmospheric science provides another important motivation. The dispersion of airborne pollutants is classically described by advection--diffusion models \cite{Stockie2011Mathematics}, and fractional derivative models have been introduced to reproduce the heavy-tailed concentration profiles observed in field measurements \cite{Goulart2017Fractional}. Since the intensity of atmospheric turbulence varies with time, the corresponding transport exponent is naturally time-dependent, which leads to variable-exponent diffusion models of the type studied in this paper.

Systematic and in-depth progress has been achieved in both forward and inverse problems for constant-order 
time-fractional diffusion equations. Sakamoto and Yamamoto \cite{Sakamoto2011Initial} established the 
well-posedness theory of 
fractional diffusion-wave equations on bounded domains by means of the Laplace transform and 
Mittag-Leffler functions. This theory was subsequently extended to multi-term time-fractional 
diffusion equations \cite{Li2015MultiTerm} and to a unified framework for weak and strong solutions of non-homogeneous 
initial boundary value problems \cite{Kian2021Wellposedness}. 
Luchko \cite{Luchko2012Initial,Luchko2016General,Luchko2010Some} and Zacher \cite{Zacher2012Global,Zacher2013DeGiorgi} 
further conducted in-depth 
investigations on forward problems for time-fractional diffusion equations from the perspectives of 
well-posedness and regularity, respectively, covering fundamental theoretical results such as existence 
and uniqueness of solutions, Hölder regularity, and long-time asymptotic behavior. 
On the inverse problem side, Jin and Rundell \cite{Jin2015Tutorial} provided a systematic 
survey of various inverse problems arising in anomalous diffusion processes, 
and the later monograph \cite{Kaltenbacher2023Inverse} by Kaltenbacher and Rundell offered a 
comprehensive framework ranging from theoretical foundations to 
numerical methods for inverse problems of fractional partial differential equations.
Liu, Rundell, and Yamamoto \cite{Liu2016Strong} established 
a uniqueness result for time-dependent source terms using the strong maximum principle for fractional 
diffusion equations, without requiring observation data at a specified terminal time. In the context of 
coefficient identification, Li and Yamamoto \cite{Li2019InverseCoefficients} systematically summarized related results on inverse 
coefficient problems for fractional partial differential equations in their review chapter, 
while Li, Liu, and Yamamoto \cite{Li2020InverseSource} made significant contributions 
to the determination of source terms for time-fractional diffusion equations.

The mathematical study of variable-exponent fractional models is still in a developmental stage. 
For the case of space-dependent variable-exponent, several mathematical and numerical 
results are available \cite{Kian2018VariableOrder,Le2021AlphaRobust,Zhao2019MultiDomain}. 
For time-dependent variable-exponent, Umarov and Steinberg \cite{Umarov2009VariableOrder} addressed 
the special case of piecewise constant variable-exponent and provided the solution to 
the Cauchy problem by means of explicit representation formulas. 
Later, Cuesta \cite{Cuesta2021SubDiffusion} studied the well-posedness of the initial value problem 
for sub-diffusion with time-dependent variable-exponent. 
Kian \cite{Kian2025PiecewiseConstant} further extended the analysis to time-fractional partial differential 
equations with time-dependent piecewise constant-exponent and established the corresponding 
well-posedness theory. In terms of numerical methods, 
Garrappa and Giusti \cite{2023Computational} proposed an efficient computational approach for 
exponential-type variable-exponent fractional differential equations and 
systematically analyzed the convergence and stability of the numerical 
schemes. In the context of inverse problems, Jiang and Li \cite{Jiang2022Coefficient} investigated 
the inverse problem of determining a time-dependent potential in a variable-exponent 
time-fractional diffusion equation, proved uniqueness based on the coercivity of the variable-exponent fractional derivative, 
and developed an iterative thresholding algorithm for numerical reconstruction. For the case of smoothly 
time-varying fractional exponent, Zheng \cite{Zheng2025TwoMethods} systematically 
proposed two general methods the convolution 
method and the perturbation method-for handling variable-exponent fractional models, applicable to 
variable-exponent initial value problems as well as more complex two-sided space-fractional boundary 
value problems and Abel integral equations. On this basis, Li, Sun, and Zheng \cite{Li2025InverseSource} studied the 
inverse source problem for a variable-exponent sub-diffusion model using the perturbation method, 
proved uniqueness by appealing to the weak unique continuation property, and established conditional stability 
through a variational identity.
Nevertheless, the theoretical study of inverse source problems for models simultaneously 
featuring time-dependent diffusion coefficients and variable-exponent fractional derivatives 
remains quite limited to date.

Based on the above discussions, we consider the inverse source problem for the following variable-exponent subdiffusion equation:
\begin{equation}\label{eq:variable-exponent-subdiffusion}
\begin{cases}
\displaystyle \partial_t^{\alpha(t)} u(x,t) - \mathcal{A}  (t) u(x,t) = F(x,t), & (x,t)\in \Omega_T := \Omega\times(0,T),\\[6pt]
u(x,0)=u_0, & x\in\Omega,\\[4pt]
u(x,t)=0, & (x,t)\in\partial\Omega\times(0,T).
\end{cases}
\end{equation}
Here $\Omega\subset\mathbb{R}^d$ with $1\le d\le 3$ is a simply-connected bounded domain with 
the piecewise smooth boundary $\partial\Omega$ with convex corners, $0<T<\infty$ is the terminal time,
 $x:=(x_1,\dots,x_d)$ denotes the spatial variables. 
 The fractional derivative of order $0<\alpha(t)<1$ is defined in terms of the variable-exponent Abel kernel $k$
\[
 \partial_t^{\alpha(t)} u := k * \partial_t u, \qquad
k(t) := \frac{t^{-\alpha(t)}}{\Gamma(1-\alpha(t))},
\]
namely,
\[ \partial_t^{\alpha(t)} u(x,t) = \frac{1}{\Gamma(1-\alpha(t))} \int_0^t (t-\tau)^{-\alpha(t)} \, \partial_\tau  u(x,\tau) \, d\tau,\]
where $\Gamma(\cdot)$ is the Gamma function.  We define a uniform elliptic operator $\mathcal{A}(t)$ by
\begin{align*}
\mathcal{A}(t)u=-\sum_{i,j=1}^{d}\partial_i\big(a_{ij}(x,t)\,\partial_j u\big),\qquad u\in H^2(\Omega),    
\end{align*}
where $\partial_j u:=\partial_{x_j} u$, and $a_{ij}=a_{ji}$, $1\le i,j\le d$, are assumed to be sufficiently smooth (e.g. $a_{ij}\in C^2([0,T];C^1(\overline\Omega))$), and satisfy
\begin{equation*}
\sum_{i,j=1}^{d} a_{ij}(x,t)\,\xi_i\xi_j\geq \rho\sum_{j=1}^{d}\xi_j^2  ,\qquad (x,t)\in\overline\Omega\times[0,T],\ \ \xi\in\mathbb{R}^d,
\end{equation*}
where $\rho>0$ is a constant independent of $x,t,\xi$. 
\begin{problem}[\textbf{Inverse source problem}]\label{prob:inverse-source}
Assume that the source term is separable as \(F(x,t)=f(x)h(t)\), 
where \(f\) is known and \(h\in L^2(0,T)\) is unknown.
 Given the observation data
\begin{align}\label{eq:observation}
    E(t):=\int_{\Omega} u(x,t)\phi(x)\,dx, \quad t\in[0,T],
\end{align}
with \(\phi\in C_0^\infty(\omega)\) and \(\omega\subset\Omega\) an 
interior observation subdomain, determine \(h(t)\).
\end{problem}

 The main contributions of this paper are as follows.
\begin{itemize}
\item[$\bullet $]  Based on the perturbation reformulation proposed in \cite{Zheng2024}, 
    we reformulate problem \eqref{eq:variable-exponent-subdiffusion} equivalently as a constant-order 
    fractional diffusion equation of 
    order \(\alpha_0=\alpha(0)\) with an additional convolution perturbation term. By using 
    the Galerkin method, the Caputo energy inequality, Young's convolution inequality, 
    and the fractional Gronwall inequality, we establish the existence, uniqueness, 
    and energy estimates for strong and weak solutions of the forward problem.
 \item[$\bullet $] Assume that the source term is separated as \(F(x,t)=f(x)h(t)\) and that \(u_0=0\). If
\(
B_f:=\int_\Omega f(x)\phi(x)\,dx\neq 0,
\)
then the unknown temporal source factor \(h(t)\) is unique and satisfies the Lipschitz stability estimate
\(
\|h\|_{L^2(0,T)}\leq C\|\partial_t^{\alpha_0}E(t)\|_{L^2(0,T)}.
\)
\item[$\bullet $] We reformulate the inverse source problem as an optimization problem 
in the Tikhonov regularization framework, construct the corresponding adjoint 
system, and propose an efficient iterative reconstruction algorithm based on 
the conjugate gradient method. For noise level \(\delta\) and regularization 
parameter \(\lambda\sim\delta^2\), we prove that the reconstruction error satisfies
\(
\|h_\lambda-h^*\|_{L^2(0,T)}\le C\delta.
\)
Numerical experiments verify the effectiveness and stability of the 
proposed algorithm for source terms of varying regularity.
\end{itemize}

The rest of this paper is organized as follows.
 Section \ref{forward_problem} is devoted to the preliminaries and the well-posedness of 
 the forward problem, where we establish the existence, uniqueness, and energy 
 estimates for both strong and weak solutions. 
 In Section \ref{inverse_problem}, we prove the uniqueness and Lipschitz stability of the inverse source problem. 
 In Section \ref{reconstruction_algorithm}, we transform the inverse problem into an optimization problem within the Tikhonov regularization framework, propose an efficient iterative reconstruction algorithm based on the 
 conjugate gradient method, and present numerical experiments to demonstrate its effectiveness and stability.

\section{Preliminaries and well-posedness of the forward problem}\label{forward_problem}

To begin with, we fix some notation used throughout the paper. Let \(L^2(\Omega)\) be the usual Lebesgue space of 
square-integrable functions on \(\Omega\) with the inner product \((\cdot,\cdot)\) and the norm \(\|\cdot\|_{L^2(\Omega)}\).
We use \(H^m(\Omega)=W^{m,2}(\Omega)\) for the standard Sobolev spaces,
and \(H_0^1(\Omega)\) for the closure of \(C_0^\infty(\Omega)\) in
\(H^1(\Omega)\).
Let \(X\) be a Banach space. For \(1\leq p<\infty\), we use
\(L^p(0,T;X)\) to denote the space of \(X\)-valued functions \(v\)
satisfying
\[
    \|v\|_{L^p(0,T;X)}
    :=\left(\int_0^T \|v(t)\|_X^p\,dt\right)^{1/p},
\]
while
\[
    \|v\|_{L^\infty(0,T;X)}
    :=\operatorname*{ess\,sup}_{0<t<T}\|v(t)\|_X.
\]


For two time-dependent functions \(g\) and \(v\), we define their
time convolution by
\[
    (g*v)(t):=\int_0^t g(t-\tau)v(\tau)\,d\tau,
    \qquad 0<t<T.
\]
For later use in the adjoint problem, we also introduce the backward
convolution
\[
    (g\circledast v)(t)
    :=\int_t^T g(s-t)v(s)\,ds,
    \qquad 0<t<T.
\]
Next, for \(\alpha>0\), we introduce the left-sided
Riemann-Liouville fractional integral operator
\[
    (J_{0+}^\alpha g)(t)
    :=\frac{1}{\Gamma(\alpha)}
      \int_0^t (t-\tau)^{\alpha-1}g(\tau)\,d\tau,
    \qquad 0<t<T.
\]
 For
\(0<\alpha<1\), the Caputo fractional derivative is defined by
\(
   \partial_t^\alpha g
    :=J_{0+}^{1-\alpha}\partial_t g.
\)
Correspondingly, the right-sided Riemann-Liouville fractional integral
and the right-sided Caputo fractional derivative are defined by
\[
    (J_T^\alpha g)(t)
    :=\frac{1}{\Gamma(\alpha)}
      \int_t^T(\tau-t)^{\alpha-1}g(\tau)\,d\tau
\]
and
\[
    D_T^\alpha g:=-J_T^{1-\alpha}\partial_t g,
\]
respectively.

By a direct calculation (or see, e.g., \cite{Jin2021Fractional}), we have the following adjoint relation:
\begin{equation}\label{Back_inte}
 \int_{0}^{T}\left(J_{0+}^{\alpha }w (t) \right)\varphi (t)  \,dt=\int_{0}^{T}w(t)\left(J_{T}^{\alpha } \varphi(t) \right)  \,dt , \quad \forall w,\, \varphi\in L^2(0,T).   
\end{equation}




In the present article, we assume that the variable exponent \(\alpha(t)\) satisfies the following conditions:
\begin{align}\label{assumption_alpha}
 \alpha(t)\in C^1[0,T],\quad \alpha(t)\in[a,b]\subset (0,1),\quad t\in[0,T].
\end{align}

\subsection{Well-posedness of the forward problem} \label{sect2}

Following the perturbation reformulation developed in \cite{Zheng2024} and subsequently used in \cite{Li_2025}, we first establish the following equivalence, which serves as a starting point for the subsequent analysis.

\begin{lemma}\label{lem_g}
    Assumption \eqref{assumption_alpha} holds and let \(\alpha_0:=\alpha(0)\). For 
 \(
u\in W^{1,1}(0,T;L^2(\Omega))
 \cap L^p(0,T;H^2(\Omega))
\)
with
\(
 p\in\left(1,\frac{1}{1-\alpha_0}\right),
 \)
 define
 \[
 \widetilde g(t)
 =
 \int_0^t
 \frac{t^{-\alpha(r)}}{\Gamma(1-\alpha(r))}
 \left[-\alpha'(r)\ln t+
 \psi (1-\alpha(r))\alpha'(r)\right]\,dr,
 \]
where \(\psi (r)=\frac{\Gamma'(r)}{\Gamma(r)}\) is the digamma function.
 Then \(u\) solves \eqref{eq:variable-exponent-subdiffusion} if and only if it solves
\begin{equation}\label{eq1:variable-exponent1}
\begin{cases}
\displaystyle \!\partial_t^{\alpha_0} u - \mathcal{A}  (t) u + \widetilde g' * u =   F+\widetilde g(t)u_0,
& (x,t)\in\Omega_T,\\[6pt]
u(x,0)=u_0, & x\in\Omega,\\[4pt]    
u(x,t)=0, & (x,t)\in\partial\Omega\times(0,T).
\end{cases}
\end{equation}
Furthermore, the equivalence between \eqref{eq:variable-exponent-subdiffusion} and \eqref{eq1:variable-exponent1} remains 
valid in the weak sense for 
\(u\in L^\infty(0,T;L^2(\Omega))\cap L^2(0,T;H_0^1(\Omega))\) satisfying
\(\partial_t^{\alpha_0}u\in L^2(0,T;H^{-1}(\Omega))\).
\end{lemma}

\begin{lemma}\label{kernel_estimate}

Let \(\alpha \) satisfy \eqref{assumption_alpha}. 
Then there exists a constant \(C>0\) such that
 \[
 |\widetilde{g}(t)| \leq C \max\left\{t^{1-\bar{\alpha}}, t^{1-\underline{\alpha}}\right\} (1 + |\ln t|), \quad 0 < t \leq T, 
\]
 and
 \[
 |\widetilde{g}'(t)| \leq C \max\left\{t^{-\bar{\alpha}}, t^{-\underline{\alpha}}\right\} (1 + |\ln t|), \quad 0 < t \leq T.
\]
Here \(\underline{\alpha} := \min_{t \in [0,T]} \alpha(t)\) and \(\bar{\alpha} := \max_{t \in [0,T]} \alpha(t)\). 
The constant \(C\) depends only on \(T\) and \(\|\alpha'\|_{L^\infty[0,T]}\).

\end{lemma}

\begin{proof}
 The proof is analogous to that of \cite[Lemma 2.2]{Li_2025}.   
\end{proof}




Based on the above estimates for the perturbation kernel, we establish the 
well-posedness of problem \eqref{eq:variable-exponent-subdiffusion} in the following two lemmas.

\begin{lemma}\label{lem_wellposedness}
Suppose $F\in W^{1,1}(0,T;L^2(\Omega))$and $u_0\in L^2(\Omega)$. Then the initial
boundary value problem \eqref{eq:variable-exponent-subdiffusion} admits a unique solution in
\(
W^{1,p}(0,T;L^2(\Omega))\cap L^p(0,T;H^2(\Omega))\cap H^1_0(\Omega)
\)
for $1<p<\dfrac{1}{1-\alpha_0}$, and the following estimate holds:
\[
\|u\|_{W^{1,p}(0,T;L^2(\Omega))}+\|u\|_{L^p(0,T;H^2(\Omega))}\le C( \|F\|_{W^{1,1}(0,T;L^2(\Omega))}+\|u_0\|_{H^2(\Omega)}),
\]
where the constant $C_T>0$ depends on $T,\Omega,d,\alpha,p$.
\end{lemma}

\begin{lemma}[Weak well-posedness]\label{lem:weak_wellposed}
Suppose \(F\in L^\infty(0,T;L^2(\Omega))\) and  \(u_0 \in L^2(\Omega)\). Then problem \eqref{eq:variable-exponent-subdiffusion} 
admits a unique weak solution
\(
u\in L^\infty(0,T;L^2(\Omega))\cap L^2(0,T;H_0^1(\Omega)),
\)
with
\(
\partial_t^{\alpha(\cdot)} u \in L^2(0,T;H^{-1}(\Omega)).
\)
Then the following energy estimate holds:
\[
\|u\|_{L^\infty(0,T;L^2(\Omega))} + \|u\|_{L^2(0,T;H_0^1(\Omega))} 
+ \|\partial_t^{\alpha(\cdot)} u\|_{L^2(0,T;H^{-1}(\Omega))} 
\leq C \bigl(\|u_0\|_{L^2(\Omega)} + \|F\|_{L^\infty(0,T;L^2(\Omega))}\bigr),
\]
where the constant \(C >0\) depends on \(T,\Omega,\alpha\).
\end{lemma}

 Lemma \ref{lem_wellposedness} has already been established in \cite{Zheng2024}, so we omit its
proof and focus on the proof of Lemma \ref{lem:weak_wellposed}.
\begin{proof}
We first use the equivalent formulation \eqref{eq1:variable-exponent1} established above.
By Lemma \ref{kernel_estimate}, we have
\[
\widetilde g'\in L^1(0,T),
\qquad
\widetilde g\in L^2(0,T).
\]

Let $\{\varphi_k\}_{k\geq 1}$ be an $L^2(\Omega)$-orthonormal basis of $H_0^1(\Omega)$ and define
\(
V_m:=\operatorname{span}\{\varphi_1,\ldots,\varphi_m\}.
\)
Let $\mathcal{P} _m$ denote the $L^2(\Omega)$-orthogonal projection onto $V_m$.
We seek an approximate solution of the form
\[
u_m(t)=\sum_{k=1}^m c_{m,k}(t)\varphi_k,
\qquad
u_m(0)=\mathcal{P}_m u_0.
\]
The Galerkin scheme is defined as follows: find \(u_m \in V_m\) such that for all \(v \in V_m\) 
and \(t \in (0,T)\), the following equation holds:  
\begin{equation}\label{eq:lemma24-galerkin}
\begin{aligned}
&\bigl(\partial_t^{\alpha_0}u_m(t),v\bigr)_{L^2(\Omega)}+a(t;u_m(t),v)
+\bigl((\widetilde g' *u_m)(t),v\bigr)_{L^2(\Omega)}\\
=&
\bigl(F(t),v\bigr)_{L^2(\Omega)}
+\widetilde g(t)(\mathcal{P}_m u_0,v)_{L^2(\Omega)}, \quad \forall v\in V_m, \quad t\in(0,T),
\end{aligned}
\end{equation}
where
\(
a(t;w,\varphi ):=\sum_{i,j=1}^d\int_\Omega
 a_{ij}(x,t)\,\partial_j w(x,t)\,\partial_i \varphi(x)\,dx, \, w,\varphi\in H_0^1(\Omega).
\)

Taking \(v=u_m(t)\) in \eqref{eq:lemma24-galerkin} and applying the Caputo energy inequality, 
the coercivity of \(a(t;\cdot,\cdot)\), Young's convolution inequality, 
and the fractional Grönwall inequality, we readily obtain the following estimate:
\begin{equation}\label{eq:lemma24-uniform-estimate}
\sup_{0\leq t\leq T}
\|u_m(t)\|_{L^2(\Omega)}^2
+
\|u_m\|_{L^2(0,T;H_0^1(\Omega))}^2
\leq
C\left(
\|u_0\|_{L^2(\Omega)}^2
+
\|F\|_{L^\infty(0,T;L^2(\Omega))}^2
\right),
\end{equation}
where \(C>0\) is independent of \(m\).

Moreover, it follows from \eqref{eq:lemma24-galerkin} that 
\begin{equation}\label{equive_lim}
  \partial_t^{\alpha_0}u_m
=
F+\mathcal{A}(t)u_m-\widetilde g'*u_m
+\widetilde g(\cdot)\mathcal{P}_m u_0.  
\end{equation}
Since \(\mathcal{A}(t):H_0^1(\Omega)\to H^{-1}(\Omega)\) is uniformly bounded for
\(t\in[0,T]\), it follows from \eqref{eq:lemma24-uniform-estimate} and
Young's convolution inequality that
\begin{equation}\label{eq:lemma24-derivative-estimate}
\|\partial_t^{\alpha_0}u_m\|_{L^2(0,T;H^{-1}(\Omega))}
\leq
C\left(
\|u_0\|_{L^2(\Omega)}
+
\|F\|_{L^\infty(0,T;L^2(\Omega))}
\right).
\end{equation}
The estimates
\eqref{eq:lemma24-uniform-estimate} and
\eqref{eq:lemma24-derivative-estimate} imply that, after passing to a
subsequence,
\[
u_m\rightharpoonup u
\quad\text{weakly in }L^2(0,T;H_0^1(\Omega)),
\]
\[
\partial_t^{\alpha_0}u_m
\rightharpoonup \zeta
\quad\text{weakly in }L^2(0,T;H^{-1}(\Omega)).
\]
Since \(H_0^1(\Omega)\) is compactly embedded in \(L^2(\Omega)\), and \(L^2(\Omega)\) 
is continuously embedded in \(H^{-1}(\Omega)\), the compactness theorem \cite{Leoni2023} for fractional energy spaces gives
\[
u_m\to u
\quad\text{strongly in }L^2(0,T;L^2(\Omega)).
\]
Consequently,
\[
\|\widetilde g'*(u_m-u)\|_{L^2(0,T;L^2(\Omega))}
\leq
\|\widetilde g'\|_{L^1(0,T)}
\|u_m-u\|_{L^2(0,T;L^2(\Omega))}
\longrightarrow0.
\]
Therefore, passing to the limit in \eqref{equive_lim}, we obtain
\[
\partial_t^{\alpha_0}u-\mathcal{A}(t)u+\widetilde g'*u
=
F+\widetilde g(t)u_0
\]
in $L^2(0,T;H^{-1}(\Omega))$. The trace theorem for the fractional energy space gives
\[
u\in L^\infty(0,T;L^2(\Omega)),
\qquad
u(0)=u_0.
\]
Thus $u$ is a weak solution of the problem \eqref{eq1:variable-exponent1}.

Let $u_1$ and $u_2$ be two weak solutions and set
\[
w:=u_1-u_2.
\]
Then
\[
\partial_t^{\alpha_0}w-\mathcal{A}(t)w+\widetilde g'*w=0,
\qquad
w(0)=0.
\]
Repeating the estimate above with $F=0$ and \(u_0=0\), we obtain
\[
\sup_{0\leq t\leq T}
\|w(t)\|_{L^2(\Omega)}^2
+
\|w\|_{L^2(0,T;H_0^1(\Omega))}^2
=0.
\]
Hence the weak solution is unique.
This completes the proof.
\end{proof}


\section{Uniqueness and Stability of the inverse source problem}\label{inverse_problem}

The assumption \(F\in L^\infty(0,T;L^2(\Omega))\) in Section \ref{sect2} is a sufficient condition for general sources. 
For the separated source \(F=f(x)h(t)\), the regularity can be relaxed to \(h\in L^2(0,T)\), 
thanks to the smoothing effects of the spatial factor and the fractional integral. 
Consequently, all subsequent theoretical analyses and numerical reconstructions are carried out under the 
setting \(F=f(x)h(t)\) with \(h\in L^2(0,T)\).

\begin{lemma}\label{lem:L2estimate}
 Let \(u\) be the unique weak solution to problem \eqref{eq:variable-exponent-subdiffusion} 
 with \(u_0 = 0\) and  $F(x,t)=f(x)h(t)$, where 
 \(f\in L^2(\Omega)\) and \(h\in L^2(0,T)\).
 Then the solution satisfies the following  estimate:
\[
\| u(\cdot,t)\|_{L^2(\Omega)}^2 \leq C J_{0_+}^{\alpha_0}|h(t)|^2, \quad  \text{for a.e. } t\in(0,T),
\]
where the constant $C>0$ depends on $T,\Omega,\alpha$, but is independent of $u$ and $h$.
\end{lemma}


\begin{proof}
Firstly, by noting that \eqref{eq:variable-exponent-subdiffusion} is equivalent to \eqref{eq1:variable-exponent1},
we multiply $u$ on both sides of the first equation of \eqref{eq1:variable-exponent1} with $u_0=0$ and integrate over $\Omega$.
Integration by parts yields
\[
\int_{\Omega} u\,\partial_t^{\alpha_0}u\,dx
+ \|\nabla u\|_{L^2(\Omega)}^2
+ \int_{\Omega} u(t)\,(\widetilde g' * u)(t)\,dx
= \int_{\Omega}  F u(t)\,dx,\qquad 0<t<T.
\]
Next, we estimate the left-hand side of the above equation from below. 
By the standard inequality for the Caputo derivative,
\[
2\,\bigl(\partial_t^{\alpha_0}u(t),u(t)\bigr)_{L^2(\Omega)}
\geq \,\partial_t^{\alpha_0}
\|u(t)\|_{L^2(\Omega)}^2,
\]
together with the Poincaré inequality, we arrive at

\begin{align}\label{J0+}
\frac{1}{2}\partial_t^{\alpha_0}\|u\|_{L^2(\Omega)}^2 + \lambda\|u\|_{H^1(\Omega)}^2
&\le \left|\int_{\Omega} u(t)(\widetilde g' * u)(t)\,dx\right| + \left|\int_{\Omega}F u(t)\,dx\right|,
\end{align}
 where \(\lambda>0\) depends only on \(\Omega\), \(d\), and \(\rho\).
Applying the fractional integral operator \(J_{0+}^{\alpha_0}\)
to both sides of \eqref{J0+}, and using
\[
J_{0+}^{\alpha_0}
\left(\partial_t^{\alpha_0}\varphi\right)(t)
=\varphi(t)-\varphi(0),
\]
together with the assumption \(u(0)=0\), we obtain
\begin{align}\label{eq:2.3}
\frac{1}{2}\|u\|_{L^2(\Omega)}^2 + \lambda J_{0+}^{\alpha_0}\|u\|_{H^1(\Omega)}^2
&\le  J_{0+}^{\alpha_0}\Big|\int_{\Omega} u(t)(\widetilde g' * u)(t)\,dx\Big|
+ J_{0+}^{\alpha_0}\Big|\int_{\Omega}F u(t)\,dx\Big|\notag\\
&=: RHS1 + RHS2,\qquad 0<t<T.
\end{align}

We estimate \(RHS1\) and \(RHS2\) separately.
For \(RHS2\),  the Cauchy-Schwarz and Young inequalities, together with the representation
\(F(x,t)=f(x)h(t)\), yield
\begin{align}\label{RHS2}
RHS2 &\le C_\varepsilon J_{0+}^{\alpha_0}\|F\|_{L^2(\Omega)}^2 + \varepsilon J_{0+}^{\alpha_0}\|u\|_{L^2(\Omega)}^2\notag\\
&\le C_\varepsilon J_{0+}^{\alpha_0}\left\lvert h(t)\right\rvert  ^2 + \varepsilon J_{0+}^{\alpha_0}\|u\|_{H^1(\Omega)}^2.
\end{align}
By choosing \(\varepsilon>0\) sufficiently small, we can absorb the last term on the right-hand side of (\ref{RHS2}) into the left-hand side of \eqref{eq:2.3}.
For $RHS1$, by the Young's inequality and Cauchy-Schwarz inequality, we have
\begin{align}\label{RHS1}
RHS1 &= J_{0+}^{\alpha_0}\Big|\int_{\Omega} u(t)(\widetilde g' * u)(t)\,dx\Big|\notag\\
&\le J_{0+}^{\alpha_0}\|u(t)\|_{L^2(\Omega)}\|(\widetilde g' * u)(t)\|_{L^2(\Omega)}\notag\\
&\le \varepsilon J_{0+}^{\alpha_0}\|u(t)\|_{L^2(\Omega)}^2
+ C_\varepsilon J_{0+}^{\alpha_0}\|    (\widetilde g' * u)(t)\|_{L^2(\Omega)}^2.
\end{align}
By letting $\varepsilon$ be sufficiently small  so that 
$\varepsilon J_{0+}^{\alpha_0}\|u(t)\|_{L^2(\Omega)}^2$ can be absorbed by $ \lambda J_{0+}^{\alpha_0}\|u\|_{H^1(\Omega)}^2$.
It remains to estimate the second term on the right-hand side of \eqref{RHS1}.
For this purpose, we use the kernel estimate established in 
Lemma \ref{kernel_estimate}:
\[
|\widetilde g'(s)|
\leq C\max\{s^{-\overline{\alpha}},
s^{-\underline{\alpha}}\}
\bigl(1+|\ln s|\bigr),
\qquad s>0.
\]

By the Young's inequality, we have
\begin{align*} 
 \|(\widetilde g' * u)(t)\|_{L^2(\Omega)}^2 &= \Big\|\int_0^t \widetilde g'(t-\tau)u(\tau)\,d\tau\Big\|_{L^2(\Omega)}^2\notag\\
&\le \int_0^t |\widetilde g'(t-\tau)|(1+\left\lvert \ln(t-\tau)\right\rvert )\,d\tau \int_0^t \frac{|\widetilde g'(t-\tau)|}{(1+\left\lvert \ln(t-\tau)\right\rvert )}\,\|u(\tau)\|_{L^2(\Omega)}^2\,d\tau\notag\\
&\le C \int_0^t \max\{(t-\tau)^{-\overline{\alpha}},(t-\tau)^{-\underline\alpha}\}\|u(\tau)\|_{L^2(\Omega)}^2\,d\tau.\notag
\end{align*} 
Since \(\underline\alpha\leq\overline\alpha\) and \(0<t\leq T\), it follows that
\[
(t-\tau)^{-\underline\alpha}
\leq T^{\overline\alpha-\underline\alpha}
(t-\tau)^{-\overline\alpha},
\qquad 0<\tau<t.
\]
Consequently,
\[
\max\{(t-\tau)^{-\overline\alpha},
(t-\tau)^{-\underline\alpha}\}
\leq C_T(t-\tau)^{-\overline\alpha},
\]
which implies that 
\begin{align}\label{eq_convo}
     \|    (\widetilde g' * u)(t)\|_{L^2(\Omega)}^2 
  &\le C \int_0^t (t-\tau)^{-\overline{\alpha}}\|u(\tau)\|_{L^2(\Omega)}^2\,d\tau.
\end{align} 
Using the fractional integral operator \(J_{0+}^{\alpha_0}\) to
 \eqref{eq_convo} and using Fubini's theorem gives

\begin{align}\label{J_be}
J_{0+}^{\alpha_0}\|    (\widetilde g' * u)(t)\|_{L^2(\Omega)}^2 
&\le C \int_0^t \int_{\tau}^{t} (t-s)^{\alpha_0-1}(s-\tau)^{-\overline{\alpha}}  ds\, \|u(\tau)\|_{L^2(\Omega)}^2 d\tau\notag \\
&\le C \int_0^t \int_{0}^{t-\tau} (t-\tau-s)^{\alpha_0-1}s^{-\overline\alpha} ds\, \|u(\tau)\|_{L^2(\Omega)}^2\, d\tau .
\end{align} 
Since \(0<\alpha_0<1\) and \(0<\overline\alpha<1\), the inner integral
can be evaluated using the Beta-function identity:
\[
\begin{aligned}
\int_0^{t-\tau}
(t-\tau-s)^{\alpha_0-1}s^{-\overline\alpha}\,\mathrm ds
&=(t-\tau)^{\alpha_0-\overline\alpha}
\int_0^1(1-r)^{\alpha_0-1}r^{-\overline\alpha}\,\mathrm dr\\
&=B(\alpha_0,1-\overline\alpha)
(t-\tau)^{\alpha_0-\overline\alpha}.
\end{aligned}
\]
Since \(B(\alpha_0,1-\overline\alpha)\) is finite, it can be
absorbed into the generic constant \(C\). Therefore, \eqref{J_be} yields
\begin{equation}\label{bet_function}
J_{0+}^{\alpha_0}
\|(\widetilde g' *u)(t)\|_{L^2(\Omega)}^2
\leq C\int_0^t
(t-\tau)^{\alpha_0-\overline\alpha}
\|u(\tau)\|_{L^2(\Omega)}^2\,\mathrm d\tau.
\end{equation}
Combining \eqref{eq:2.3}-(\ref{RHS1}) with the estimate in \eqref{bet_function}, and choosing \(\varepsilon>0\) sufficiently small, leads to

\begin{align}\label{2_eq}
\|u(t)\|_{L^2(\Omega)}^2 
&\le  C J_{0+}^{\alpha_0}\left\lvert h(t)\right\rvert  ^2 + C \int_0^t (t-\tau)^{\alpha_0-\bar\alpha}\|u(\tau)\|_{L^2(\Omega)}^2 \, d\tau.
\end{align}

Applying the generalized fractional Grönwall  inequality to \eqref{2_eq}, we arrive at
\begin{align*}
&\|u(t)\|_{L^2(\Omega)}^2 
\le  C J_{0+}^{\alpha_0}\left\lvert h(t)\right\rvert  ^2 \notag\\
&+ C \int_0^t (t-\tau)^{\alpha_0-\bar\alpha}E_{\alpha_0-\bar{\alpha}+1,\alpha_0-\bar\alpha+1}(C\Gamma(\alpha_0-\bar\alpha+1)(t-\tau)^{\alpha_0-\bar\alpha+1})J_{0_+}^{\alpha_0}\left\lvert h (\tau)\right\rvert ^2 \, d\tau\\
&\le  C J_{0+}^{\alpha_0}\left\lvert h(t)\right\rvert  ^2 + C \int_0^t (t-\tau)^{\alpha_0-\bar\alpha}J_{0_+}^{\alpha_0}\left\lvert h (\tau)\right\rvert ^2 \, d\tau.
\end{align*}
Since \(\alpha_0-\bar\alpha>-1\), the Mittag-Leffler function in the above integral is bounded for \(0\le t-\tau\le T\). 
Hence, By the convolution formula for fractional integrals (see, e.g.,\cite[lemma 2.1 ]{hu2024strong}) and the semigroup property of fractional integrals, we have
\begin{align*}
   \|u(t)\|_{L^2(\Omega)}^2 
   &\le  C J_{0+}^{\alpha_0}\left\lvert h(t)\right\rvert  ^2 + C J_{0+}^{\alpha_0}\left(\int_0^t (t-\tau)^{\alpha_0-\bar\alpha}\left\lvert h (\tau)\right\rvert ^2 \, d\tau\right)\notag\\
   &\le  C J_{0+}^{\alpha_0}\left\lvert h(t)\right\rvert  ^2 + C J_{0+}^{2\alpha_0-\bar\alpha+1}\left\lvert h (t)\right\rvert ^2\notag\\
   &\le C  J_{0+}^{\alpha_0}\left\lvert h(t)\right\rvert  ^2.
\end{align*}
This completes the proof.
\end{proof}

\begin{theorem}\label{stability}
Let \(u_0=0\) and \(F(x,t)=f(x)h(t)\) with \(h\in L^2(0,T)\). 
Assume further that \(B_f:=\int_\Omega f(x)\phi(x)dx\neq0\). Then there exists a constant \(C>0\), independent of \(h\), such that
\[\|h\|_{L^2(0,T)} \le C \left\|\partial_t^{\alpha_0} E (t)\right\|_{L^2(0,T)}.\]
\end{theorem}
\begin{proof}
We multiply $\phi  (x)$ on both sides of the first equation of \eqref{eq1:variable-exponent1} with $u_0=0$ and integrate over $\Omega$. 
Integration by parts yields
\[\int_{\Omega} \partial_t^{\alpha_0}u(x,t)\,\phi(x)\,dx- \int_{\Omega} \mathcal{A} (t)u(x,t)\phi(x)\,dx
= h(t)\int_{\Omega}  f(x) \phi(x)\,dx-\int_{\Omega} \phi(x)(\widetilde g' * u)(x,t)\,dx.
\]
Setting $B_f=\int_{\Omega} f(x) \phi(x)\,dx $, we have
\[
h(t) = B_f^{-1}\left[\int_{\Omega} \partial_t^{\alpha_0}u(x,t)\,\phi(x)\,dx- \int_{\Omega} \mathcal{A} (t)u(x,t)\phi(x)\,dx + \int_{\Omega} \phi(x)(\widetilde g' * u)(x,t)\,dx\right]  .
\]
Squaring both sides of the above inequality and using the observation data \eqref{eq:observation}, we obtain
\begin{align*}
    \left\lvert h(t)\right\rvert^2 \leq \frac{1}{\left\lVert B_f\right\rVert^2 }\left[(\partial_t^{\alpha_0} E (t))^2 + \left(\int_{\Omega} \mathcal{A} (t)u(x,t)\phi(x)\,dx\right)  ^2 + \left(\int_{\Omega} \phi(x)(\widetilde g' * u)(x,t)\,dx\right) ^2\right] .
\end{align*}
Now we estimate the right-hand side of the above inequality. By the integration by parts formula 
and the fact that $\phi\in C_{0}^\infty(\Omega)$, we see that
\begin{align*}
   & \int_{\Omega} \mathcal{A} (t)u(x,t)\phi(x)\,dx = -\int_{\Omega}^{}\sum_{i,j=1}^{d}\partial_i\big(a_{ij}(x,t)\,\partial_j u\big)\phi\,dx \\
    &= -\int_{\Omega}^{}u\left(\sum_{i,j=1}^{d} \,\partial_j \left(a_{ij}(x,t)\,\partial_i\phi\right)\right)  \,dx \notag\\
    &=\int_{\Omega}^{} u \, \mathcal{A}(t)\phi \,dx .
\end{align*}
Combining the Cauchy-Schwarz inequality with Lemma \ref{lem:L2estimate}, we  further have
\begin{align*}
    \left(\int_{\Omega} \mathcal{A} (t)u(x,t)\phi(x)\,dx\right)  ^2 \leq \left\lVert \mathcal{A}(t)\phi \right\rVert_{L^2(\Omega)}^2 \left\lVert u(\cdot,t) \right\rVert_{L^2(\Omega)}^2 \leq C  J_{0+}^{\alpha_0}\left\lvert h(t)\right\rvert  ^2.
\end{align*}
Similarly, we can easily prove that
\begin{align}\label{second}
   & \left(\int_{\Omega} \phi(x)(\widetilde g' * u)(x,t)\,dx\right) ^2 \leq \left\lVert \phi \right\rVert_{L^2(\Omega)}^2 \left\lVert (\widetilde g' * u)(\cdot,t) \right\rVert_{L^2(\Omega)}^2 \notag\\
   &\le C \int_0^t (t-\tau)^{-\overline{\alpha}}\|u(\tau)\|_{L^2(\Omega)}^2\,d\tau \leq C \int_0^t (t-\tau)^{-\overline{\alpha}}J_{0+}^{\alpha_0}\left\lvert h (\tau)\right\rvert ^2 \, d\tau\notag\\
    &\leq C \int_0^t (t-\tau)^{\alpha_0-\overline{\alpha}}|h(\tau)|^2\,d\tau.
\end{align}
\eqref{second}  together with Lemma \ref{lem_wellposedness}, we obtain
\begin{align*}
    \left\lvert h(t)\right\rvert^2  &\leq  \frac{C}{\left\lVert B_f\right\rVert^2 }\left[(\partial_t^{\alpha_0} E (t))^2  + \int_0^t (t-\tau)^{\alpha_0-1}(1+(t-\tau)^{1-\overline{\alpha}}) |h(\tau)|^2\,d\tau\right] \\
    &\leq  \frac{C}{\left\lVert B_f\right\rVert^2 }\left[(\partial_t^{\alpha_0} E (t))^2  + (1+T)\int_0^t (t-\tau)^{\alpha_0-1} |h(\tau)|^2\,d\tau\right] \\
    &\leq  \frac{C}{\left\lVert B_f\right\rVert^2 }\left[(\partial_t^{\alpha_0} E (t))^2  + \int_0^t (t-\tau)^{\alpha_0-1} |h(\tau)|^2\,d\tau\right] .
\end{align*}
By Grönwall's inequality, we have 
\[
|h(t)|^2 \le a(t) + b \int_0^t (t-\tau)^{\alpha_0-1} E_{\alpha_0,\alpha_0}\bigl(b (t-\tau)^{\alpha_0}\bigr) a(\tau)\,d\tau,
\]
where $a(t) = C \left(\frac{\partial_t^{\alpha_0} E (t)}{\left\lVert B_f\right\rVert }\right)^2 $ and $b = C \frac{\Gamma(\alpha_0)}{\left\lVert B_f\right\rVert^2 } $. 

By the boundedness of the Mittag-Leffler function \(E_{\alpha_0,\alpha_0}(b(t-\tau)^{\alpha_0})\) on \(0 \le t-\tau \le T\), we get

\begin{align}\label{h-estimate}
|h(t)|^2 &\le a(t) + C \int_0^t (t-\tau)^{\alpha_0-1} a(\tau)\,d\tau.
\end{align}
Integrate both sides of \eqref{h-estimate} over \([0,T]\) and swap the order of the double integral on the right:
\[
\int_0^T \int_0^t (t-\tau)^{\alpha_0-1} a(\tau)\,d\tau\,dt \le \frac{T^{\alpha_0}}{\alpha_0} \int_0^T a(\tau) \,d\tau.
\]
This yields
\[
\int_0^T |h(t)|^2 dt \le (1+ C \frac{T^{\alpha_0}}{\alpha_0})\int_0^T a(t) \,dt.
\]
Consequently, 
\[\|h\|_{L^2(0,T)} \le C \left\|\partial_t^{\alpha_0} E (t)\right\|_{L^2(0,T)}.\]
This completes the proof.
\end{proof}

\begin{remark}
    The uniqueness of the inverse problem follows immediately from Theorem \ref{stability} by taking \(E(t)=0\). 
    Indeed, the stability estimate yields \(\|h\|_{L^2(0,T)}\le C\|\partial_t^{\alpha_0}E\|_{L^2(0,T)}=0\), which implies \(h=0\) a.e. in \((0,T)\).
\end{remark}

\section{ Numerical algorithm and experiments}\label{reconstruction_algorithm}
We consider Tikhonov regularization for the inverse problem of determining the source term $h(t)$ 
in \eqref{eq:variable-exponent-subdiffusion} 
from the observation data \eqref{eq:observation}.

Now let $\mathcal{K} $ be the solution operator that maps the source term $h(\cdot)$ to the observation data $E(\cdot)$, which is defined by
\[  \mathcal{K} h = \int_{\Omega} u(x,t)\phi(x)\,dx,\quad \phi\in C_{0}^{\infty}(\omega ),\]
where $u$ is the solution to \eqref{eq:variable-exponent-subdiffusion} with the source term $F(x,t)=f(x)h(t)$.

By the well-posedness result in Lemma \ref{lem:weak_wellposed}, the operator \(\mathcal{K}\) is linear and bounded from \(L^2(0,T)\) to \(L^2(0,T)\), 
and Theorem \ref{stability} further guarantees that \(\mathcal{K}\) is  one-to-one. Moreover, let \(h^*\in L^2(0,T)\) be the unknown source term, 
and \(E^* = \mathcal{K} h^*\) be the corresponding exact observation data. 
In practice, we only have access to noisy observation data \(E^\delta\) satisfying \(\|\partial_t^{\alpha_0}\left(E^\delta - E^*\right) \|_{L^2(0,T)} \le \delta\), 
where \(\delta > 0\) is the noise level.

The Tikhonov regularization is to find $h_\lambda\in L^2(0,T)$ such that
\[h_\lambda  = \arg\min_{h\in L^2(0,T)} \left\{\|\partial_t^{\alpha_0}\left((\mathcal{K} h)(t) - E^\delta(t)\right) \|_{L^2(0,T)}^2 + \lambda \|h\|_{L^2(0,T)}^2\right\},\]
where $\lambda>0$ is the regularization parameter.

\begin{theorem}\label{thm:convergence}
 Let \(h^* \in L^2(0,T)\) be the exact source term of the problem \eqref{eq:variable-exponent-subdiffusion}, 
and let \(E^* = \mathcal{K}h^*\) be the corresponding exact observation data. Assume that the noisy data \(E^\delta \in L^2(0,T)\) satisfy
\[
\|\partial_t^{\alpha_0}\left(E^\delta - E^*\right) \|_{L^2(0,T)} \le \delta,
\]
with \(\delta > 0\). For each \(\lambda > 0\), let \(h_\lambda \in L^2(0,T)\) be the unique minimizer of the Tikhonov functional
\begin{equation}\label{tikhonov_functional}
 \mathcal{J}_\lambda(h) := \|\partial_t^{\alpha_0}(\mathcal{K}h) - \partial_t^{\alpha_0}E^\delta\|_{L^2(0,T)}^2 + \lambda \|h\|_{L^2(0,T)}^2.
\end{equation}
 If the regularization parameter \(\lambda \sim \delta^2 ( \text{as }\delta \to 0)\), then 
\[
\|h_\lambda - h^*\|_{L^2(0,T)} \leq C \delta.
\]   
\end{theorem}

\begin{proof}
First, from the minimizing property of \(h_\lambda\) and the fact that \(h^* \in L^2(0,T)\), we have
\[
\mathcal{J}_\lambda(h_\lambda) \le \mathcal{J}_\lambda(h^*).
\]
Since \(\partial_t^{\alpha_0}E^*(t) = \partial_t^{\alpha_0}(\mathcal{K}h^*(t))\), it follows that
\begin{align*}
 \|\partial_t^{\alpha_0}(\mathcal{K}h_\lambda) - \partial_t^{\alpha_0}E^\delta\|_{L^2(0,T)}^2 + \lambda \|h_\lambda\|_{L^2(0,T)}^2
\le& \|\partial_t^{\alpha_0}E^\delta - \partial_t^{\alpha_0}E^*\|_{L^2(0,T)}^2 + \lambda \|h^*\|_{L^2(0,T)}^2\\
\le &\delta^2 + \lambda \|h^*\|_{L^2(0,T)}^2.   
\end{align*}

Consequently, we obtain the following two estimates:
\begin{equation}\label{estimate1}
    \|\partial_t^{\alpha_0}(\mathcal{K}h_\lambda) - \partial_t^{\alpha_0}E^\delta\|_{L^2(0,T)}^2 \le \delta^2 + \lambda \|h^*\|_{L^2(0,T)}^2,
\end{equation}
\begin{equation}\label{estimate2}
    \lambda\|h_\lambda\|_{L^2(0,T)}^2 \le \lambda\|h^*\|_{L^2(0,T)}^2 + \delta^2. 
\end{equation}

Next, we estimate the error \(\|h_\lambda - h^*\|_{L^2(0,T)}\). By the Theorem \ref{thm:convergence} and the triangle inequality, we have
\begin{align*}
    \|h_\lambda - h^*\|_{L^2(0,T)}^2 &\leq C \|\partial_t^{\alpha_0}\mathcal{K}(h_\lambda - h^*)\|_{L^2(0,T)}^2 \\
    &\le C\left(\|\partial_t^{\alpha_0}\mathcal{K}h_\lambda - \partial_t^{\alpha_0}E^\delta\|_{L^2(0,T)}^2 + \|\partial_t^{\alpha_0}(E^\delta - E^* )\|_{L^2(0,T)}^2\right).
\end{align*}
Using \eqref{estimate1}-\eqref{estimate2} and the fact that \(\|\partial_t^{\alpha_0}E^\delta - \partial_t^{\alpha_0}E^*\|_{L^2(0,T)} \le \delta\), we obtain
\begin{align*}
\|h_\lambda - h^*\|_{L^2(0,T)}^2 &\le C\left(\delta^2 + \lambda \|h^*\|_{L^2(0,T)}^2 + \delta^2\right)\\
&\le C\left(\delta^2 + \lambda \|h^*\|_{L^2(0,T)}^2\right).
\end{align*}
By choosing \(\lambda \sim \delta^2\), we have
\[\|h_\lambda - h^*\|_{L^2(0,T)}^2 \le C \delta^2.\]
This completes the proof.
\end{proof}

To make use of efficient iterative methods, we need information about the Gâteaux derivative \( \mathcal{J}'_\lambda(h)\xi \) of the Tikhonov 
functional \( \mathcal{J}_\lambda(h) \) defined in \eqref{tikhonov_functional}. 
For any fixed direction \( \xi \in L^2(0,T) \), a direct calculation yields:
\begin{equation}\label{G_derivate}
    \mathcal{J}'_\lambda(h)\xi = 2 \int_0^T R(t) \cdot \partial_t^{\alpha_0}\!\left(\int_\Omega u(\xi) \,\phi\,dx\right)\,dt + 2 \lambda \int_0^T h(t)\xi(t)\,dt,
\end{equation}
where \( R(t) = \partial_t^{\alpha_0}(\mathcal{K}h)(t) - \partial_t^{\alpha_0}E^\delta(t) \),  
and \( u(\xi) \) is the solution to \eqref{eq:variable-exponent-subdiffusion} with the source term \( F(x,t) = f(x)\xi(t) \).  
Moreover, we have used the linearity of \(u(h)\) with respect to \(h\), which implies that  \(u'(h)\xi = u(\xi)\).

However, using \eqref{G_derivate} to evaluate \(\mathcal{J}'(h)\xi\) for every \(\xi \in L^2(0,T)\) is too costly, 
because for a given \(h\), this would require solving \eqref{eq:variable-exponent-subdiffusion} to obtain \(u(\xi)\) for each direction \(\xi\) in \(L^2(0,T)\). 
Hence, in order to reduce the computational burden, by Lemma \ref{lem_g}, 
we introduce the adjoint system corresponding to \eqref{eq:variable-exponent-subdiffusion} as
\begin{equation}\label{adjoint_system}
\begin{cases}
D_T^{\alpha_0} v - \mathcal{A}(t) v + \widetilde{g}'\circledast  v = D_T^{\alpha_0} R(t) \, \phi(x), & (x,t)\in \Omega_T, \\[6pt]
(J_T^{1-\alpha_0} v)(x,T) = 0, & x\in \Omega, \\[4pt]
v(x,t) = 0, & (x,t)\in \partial\Omega \times (0,T],
\end{cases}     
\end{equation}
where $\widetilde{g}'\circledast  v = \int_t^T \widetilde{g}'(s-t) v(s) ds$, and $J_T^{\alpha_0} $ is 
the backward Riemann-Liouville fractional integral operator defined by
$$J_T^{\alpha_0} g(t) = \frac{1}{\Gamma(\alpha_0)} \int_t^T g(\tau) (\tau-t)^{\alpha_0-1}\, d\tau,$$
and the backward Riemann-Liouville fractional derivative of order $\alpha_0$ is defined by
$$
D_T^{\alpha_0} g := -J_T^{1-\alpha_0}\frac{d}{dt}\, g(t).
$$
Then, by using the \eqref{Back_inte}, it is straightforward to show that
\[
\int_0^T R(t) \cdot \partial_t^{\alpha_0}\!\left(\int_\Omega u(\xi) \,\phi\,dx\right)\,dt
= \int_0^T \int_\Omega \left(u(\xi) \cdot \phi\right) \, D_T^{\alpha_0}R(t)\,dx\,dt.
\]
Consequently, \eqref{G_derivate} can be rewritten as
\begin{equation}\label{G_derivate1}
    \mathcal{J}'_\lambda(h)\xi = 2 \int_0^T \int_\Omega u_\xi \cdot \phi\, D_T^{\alpha_0}R(t)\,dx\,dt + 2 \lambda \int_0^T h(t)\xi(t)\,dt.
\end{equation}

We now multiply both sides of \eqref{eq1:variable-exponent1} by \(v\) and integrate over \(\Omega_T\), which yields
\begin{equation}\label{integr}
 \int_0^T\int_\Omega \big(\partial_t^{\alpha_0}u_\xi - \mathcal{A}u_\xi + \widetilde{g}'*u_\xi\big) v\,dxdt = \int_0^T\int_\Omega f(x)\xi(t) v\,dxdt,   
\end{equation}
where \(u(\xi)\) is the solution to \eqref{eq1:variable-exponent1} with the source term \(F(x,t) = f(x)\xi(t)\). 
Then, applying integration by parts to \eqref{integr}, we obtain
\begin{equation*}
\int_0^T\int_\Omega u_\xi \left(D_T^{\alpha_0}v - \mathcal{A}v + \widetilde{g}'\circledast  v\right) dxdt= \int_0^T \xi(t) \int_\Omega f(x)v(x,t)\,dxdt.
\end{equation*}
Using the adjoint system \eqref{adjoint_system}, the left‑hand side of the above equation equals 
\[
  \int_0^T\int_\Omega u_\xi \phi\, D_T^{\alpha_0}R(t)\,dx\,dt = \int_0^T \xi(t) \int_\Omega f(x)v(x,t)\,dx\,dt.
\]
Substituting this expression into \eqref{G_derivate1}, then the first variation 
of \(\mathcal J_\lambda\) in an arbitrary direction \(\xi\) is given by
\[
\mathcal{J_\lambda}'(h)\xi = 2 \int_0^T \xi(t) \left( \int_\Omega f(x)v (x,t)\,dx +  \lambda h(t) \right) dt.
\]
Consequently, the \(L^2(0,T)\) gradient of the functional is
\[
\nabla\mathcal J_\lambda(h)(t)
=
2\left[
\int_\Omega f(x)v(x,t)\,dx
+
\lambda h(t)
\right].
\]

The numerical reconstruction is performed using the standard 
Tikhonov regularization and the conjugate gradient (CG) 
method (see, e.g., \cite{LI2026multifactor}). The detailed 
iteration process is described as follows.

\begin{algorithm}
Choose an initial guess \(h^0 \in L^2(0,T)\). The iteration process is then given by
 \[
 h^{k+1}=h^k + r^k d^k,\quad k=0,1,\dots
 \]
 where the first search direction is defined by
 \[
 d^0 = -\nabla\mathcal J_\lambda(h^0),
 \]
 and, for \(k\ge 1\),
 \[
 d^k = -\nabla\mathcal J_\lambda(h^k) + s^k d^{k-1}.
 \]
 The Fletcher-Reeves coefficient is defined as
 \[
 s^k = \frac{\|\nabla\mathcal J_\lambda(h^k)\|_{L^2(0,T)}^2}{\|\nabla\mathcal J_\lambda(h^{k-1})\|_{L^2(0,T)}^2},
 \]
 and the step length \(r^k\) is selected by minimizing the functional along the search direction:
 \[
 r^k = \arg\min_{r  \ge 0} \mathcal J_\lambda(h^k + r d^k).
 \]
\end{algorithm}

\begin{remark}\label{rem:step-length}
    Since the forward operator \(\mathcal K\) is linear with respect to the temporal source \(h\), the functional \(\mathcal J_\lambda(h^k+r d^k)\) is a strictly convex quadratic function of \(r\). Define
 \[
 q^k=\partial_t^{\alpha_0}(\mathcal K h^k-E^\delta),
 \qquad
 z^k=\partial_t^{\alpha_0}\mathcal Kd^k.
 \]
 Then the unconstrained minimizer is given by
 \[
 \widehat r^k
 =-
 \frac{
 \bigl(q^k,z^k\bigr)_{L^2(0,T)}
 +\lambda\bigl(h^k,d^k\bigr)_{L^2(0,T)}
 }{
 \|z^k\|_{L^2(0,T)}^2
 +\lambda\|d^k\|_{L^2(0,T)}^2
 }.
 \]
 Hence, under the constraint \(r\ge0\), we take
 \[
 r^k=\max\{0,\widehat r^k\}.
 \]
 In particular, no additional iterative line-search procedure is required.
\end{remark}

\subsection{Numerical Experiments}

 In this section, we present numerical experiments to illustrate the performance of 
 the proposed reconstruction method for the inverse source problem \eqref{eq:variable-exponent-subdiffusion}.
The computational domain is the unit square \(\Omega=(0,1)^2\) over the time interval \([0,T]\) with \(T=1\). 
For spatial discretization, we employ the standard Galerkin finite element method on uniform triangular 
meshes with piecewise linear basis functions, using \(N_{xy}=20\) subdivisions in each direction. 
The variable-exponent fractional derivative is approximated by the \(L_1\) scheme on a uniform time grid with \(N_t=60\) steps.
Homogeneous Dirichlet and zero initial conditions are imposed.
At each conjugate gradient iteration, only one forward problem and one  adjoint problem need to be solved, 
which makes the numerical implementation simple and computationally inexpensive.

The spatially varying diffusion coefficient is chosen in the form 
\(\mathcal{A}(t)u=-\nabla\cdot\bigl(a(x,y,t)\nabla u\bigr)\), where \(a(x,y,t)=1+2x+2y+\frac{t}{2}\).
In the numerical tests, the source term is separated as $F(x,y,t)=f(x,y)h(t)$ with $f(x,y)=30+15\pi(x+y)$.
For the interior observation, we take the subdomain
\(\omega=\left(\frac14,\frac34\right)^2\) and use the smooth weight function
\(\phi\) supported in \(\omega\), defined by
\[
\phi(x,y)=
\begin{cases}
\displaystyle
\exp\!\left(
-\frac{1}{1-(4x-2)^2}
-\frac{1}{1-(4y-2)^2}
\right),
& (x,y)\in\omega,\\[2ex]
0, & (x,y)\notin\omega.
\end{cases}
\]

We generate noisy data in a way that is consistent with the 
theoretical setting of Theorem \ref{thm:convergence} . 
Specifically, the noisy data 
\(E^\delta\)
  are generated such that the relative noise level is
 \[ \delta_{rel}:=
 \frac{
 \left\|\partial_t^{\alpha_0}(E^\delta-E^*)\right\|_{L^2(0,T)}
 }{
 \left\|\partial_t^{\alpha_0}E^*\right\|_{L^2(0,T)}
 }
,
 \]
 where \(E^*\) and \(E^\delta\) denote the exact and noisy observations, respectively.
Unless otherwise stated, the regularization parameter is fixed as
$\lambda=5.0\times10^{-6}$.
 The reconstruction accuracy is measured by the relative $L^2(0,T)$-norm error
\[
\operatorname{error}
:=
\frac{\|h_{\rm rec}-h^*\|_{L^2(0,T)}}
{\|h^*\|_{L^2(0,T)}}.
\]

\begin{example}\label{ex1}
     Consider the problem \eqref{eq:variable-exponent-subdiffusion} 
     with the variable fractional order \(\alpha(t) = 0.5 + 0.2t\) for $t\in [0,1]$.
    we test the following cases:
\begin{itemize}
    \item[\bf {(1a):}] \(h(t)=2\sin(2\pi t)\) (left in Figure~\ref{fig:ex1});
    \item[\bf {(1b):}] \(h(t)=2\cos(3\pi t)\) (right in Figure~\ref{fig:ex1}).
\end{itemize}

  \begin{figure}[H]
       \centering
       \begin{minipage}[t]{0.48\textwidth}
       \centering
       \includegraphics[width=7.5cm]{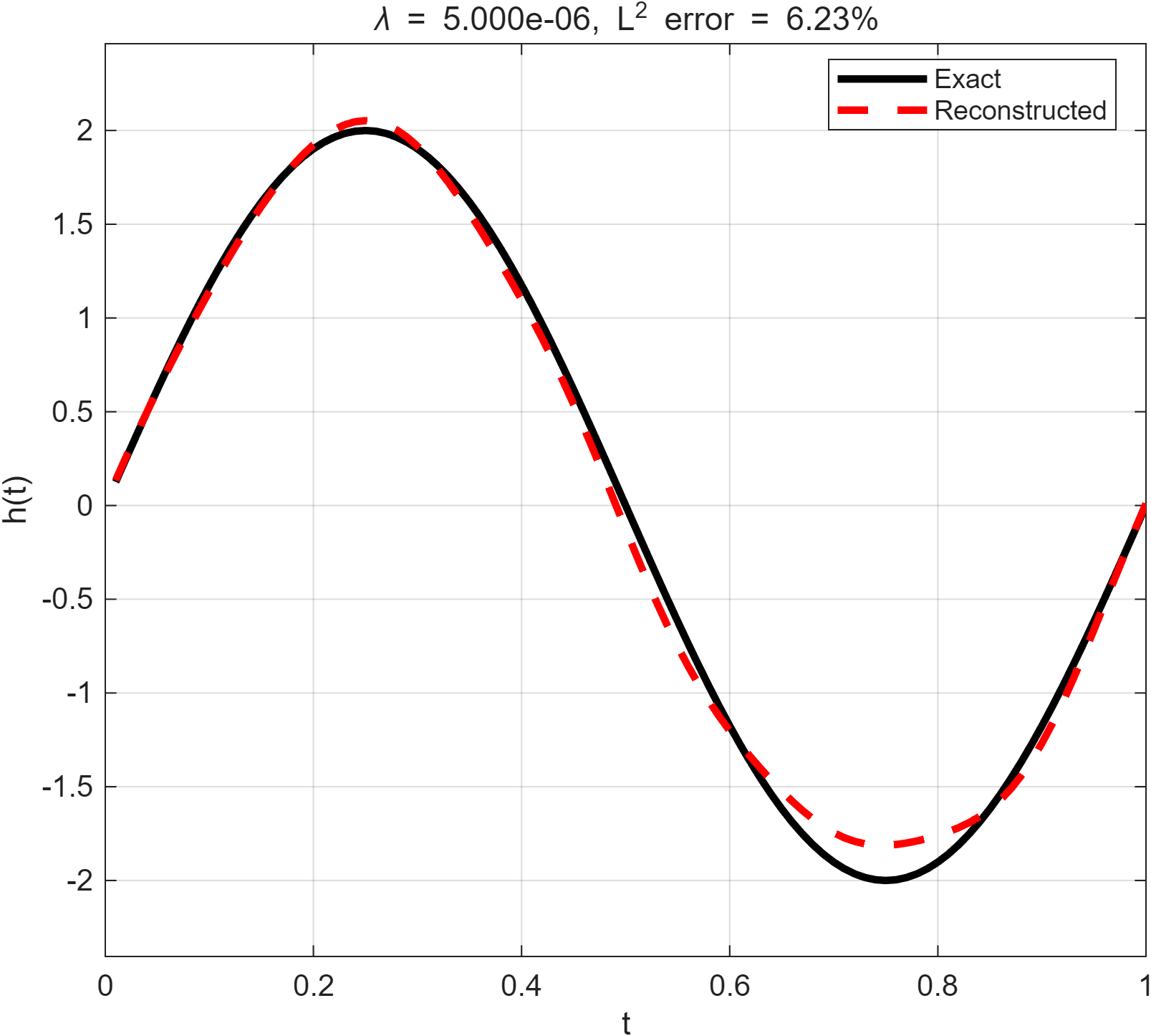}
       \end{minipage}
       \begin{minipage}[t]{0.48\textwidth}
       \centering
       \includegraphics[width=7.5cm]{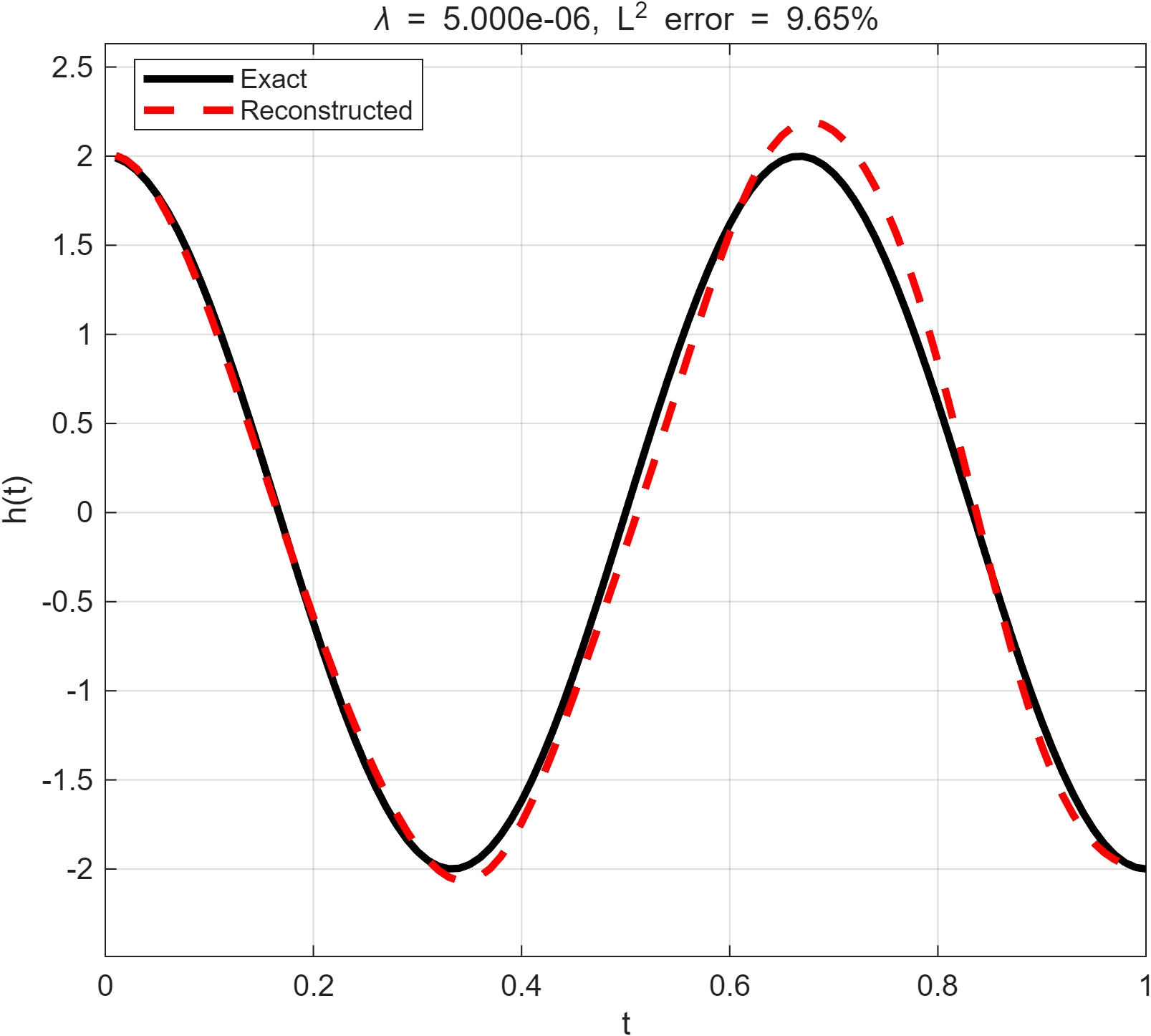}
       \end{minipage}
       \caption{Example \ref{ex1}: Numerical results for the variable fractional order \(\alpha(t) = 0.5 + 0.2t\) 
       with $\delta_{rel} = 10\%$: \((\text{left})\,error =6.23\%\), \( (\text{right})\, error=9.65\%\).
   }\label{fig:ex1}
       \end{figure}
\end{example}

\begin{table}[H]
    \centering
    \setlength{\belowcaptionskip}{8pt}   
    \setlength{\tabcolsep}{0.8cm} 
    \renewcommand{\arraystretch}{1.6} 
        \caption{Relative \(L^2\) reconstruction errors for different noise levels in Example \ref{ex1}.}\label{table1}
    \begin{tabular}{c c c}
        \Xhline{1.5pt}   
        $\delta_{rel}$ & error of $h(t) = 2\sin(2\pi t)$ & error of $h(t) = 2\cos(3\pi t)$ \\
        \hline            
        16\% & $9.84\% $ &15.32\% \\
        12\% & $7.40\% $ & 11.48\% \\
        10\% &  $6.23\% $ &9.65\% \\
        8\% & $4.96\%$ & 7.64\% \\
        4\% & $2.53\% $  & 3.80\% \\
        \Xhline{1.5pt}   
    \end{tabular}
\end{table}

 Figure \ref{fig:ex1} illustrates the exact and reconstructed temporal source functions 
 for Example \ref{ex1} with a noise level of \(\delta_{rel} = 10\%\).  
 Despite the presence of noise, the proposed method accurately captures the oscillatory profiles of both the sine and cosine sources.
Table \ref{table1} reports the relative \(L^2\) reconstruction errors for various noise levels.  
As \(\delta_{rel}\) decreases, the errors decrease steadily for both test cases, which is consistent with the theoretical 
convergence result in Theorem \ref{thm:convergence}.  In particular, the errors are quite satisfactory 
at a noise level of \(4\%\), and even for a relatively large noise level of \(\delta_{rel} = 10\%\), 
the numerical results shown in Figure \ref{fig:ex1} are still good.  
These numerical results confirm the stability and efficiency of the proposed reconstruction approach.

\begin{example}\label{ex2}
    Consider the problem \eqref{eq:variable-exponent-subdiffusion} 
    with the variable fractional order \(\alpha(t) = 0.5 + 0.2t^2\) for $t\in [0,1]$. we test the following cases:
\begin{itemize}
    \item[\bf {(2a):}] \(h(t)=2\sin(2\pi t)\) (left in Figure~\ref{fig:ex2});
    \item[\bf {(2b):}] \(h(t)=2\cos(3\pi t)\) (right in Figure~\ref{fig:ex2}).
\end{itemize}
  \begin{figure}[H]
       \centering
       \begin{minipage}[t]{0.48\textwidth}
       \centering
       \includegraphics[width=7.5cm]{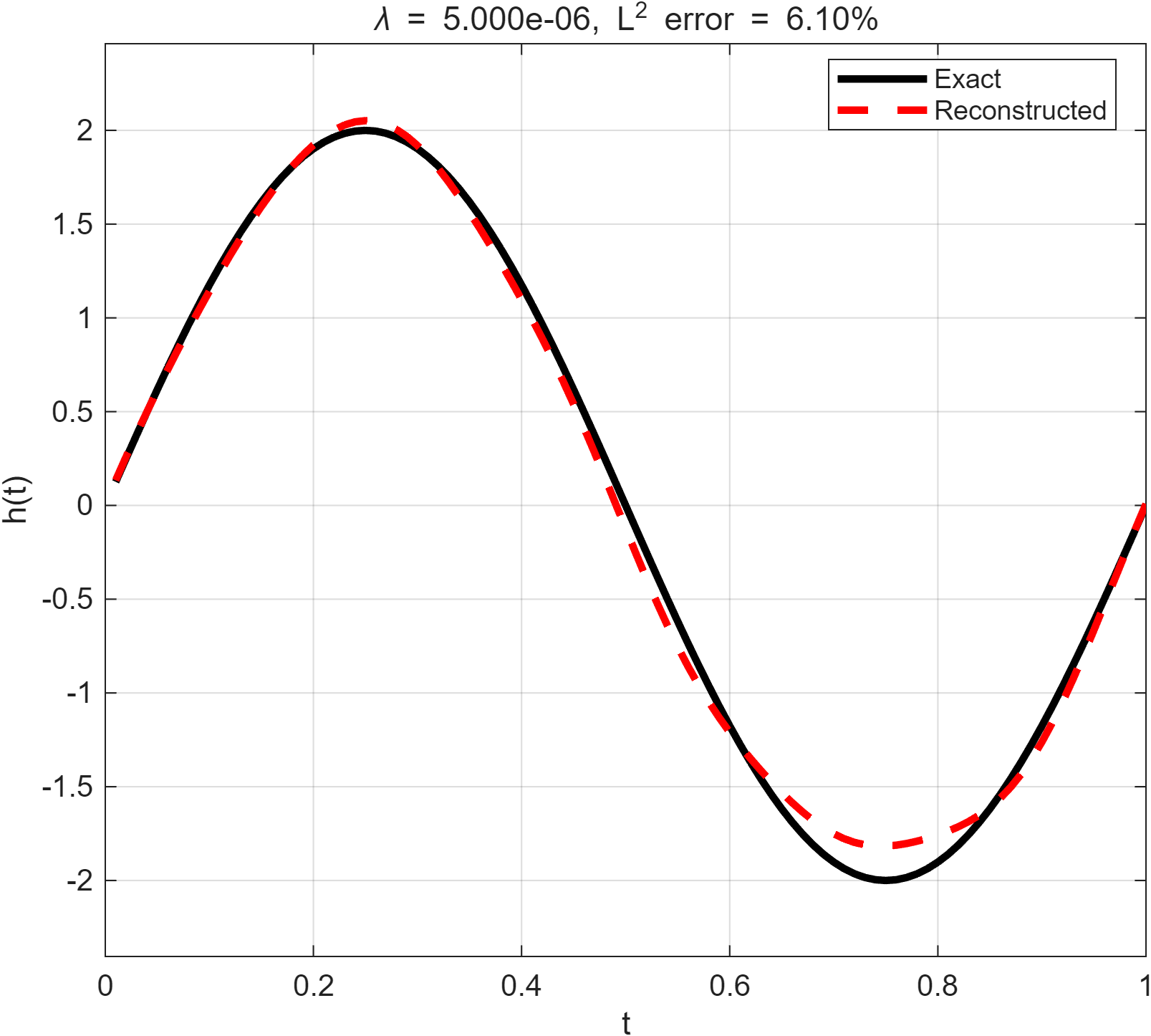}
       \end{minipage}
       \begin{minipage}[t]{0.48\textwidth}
       \centering
       \includegraphics[width=7.5cm]{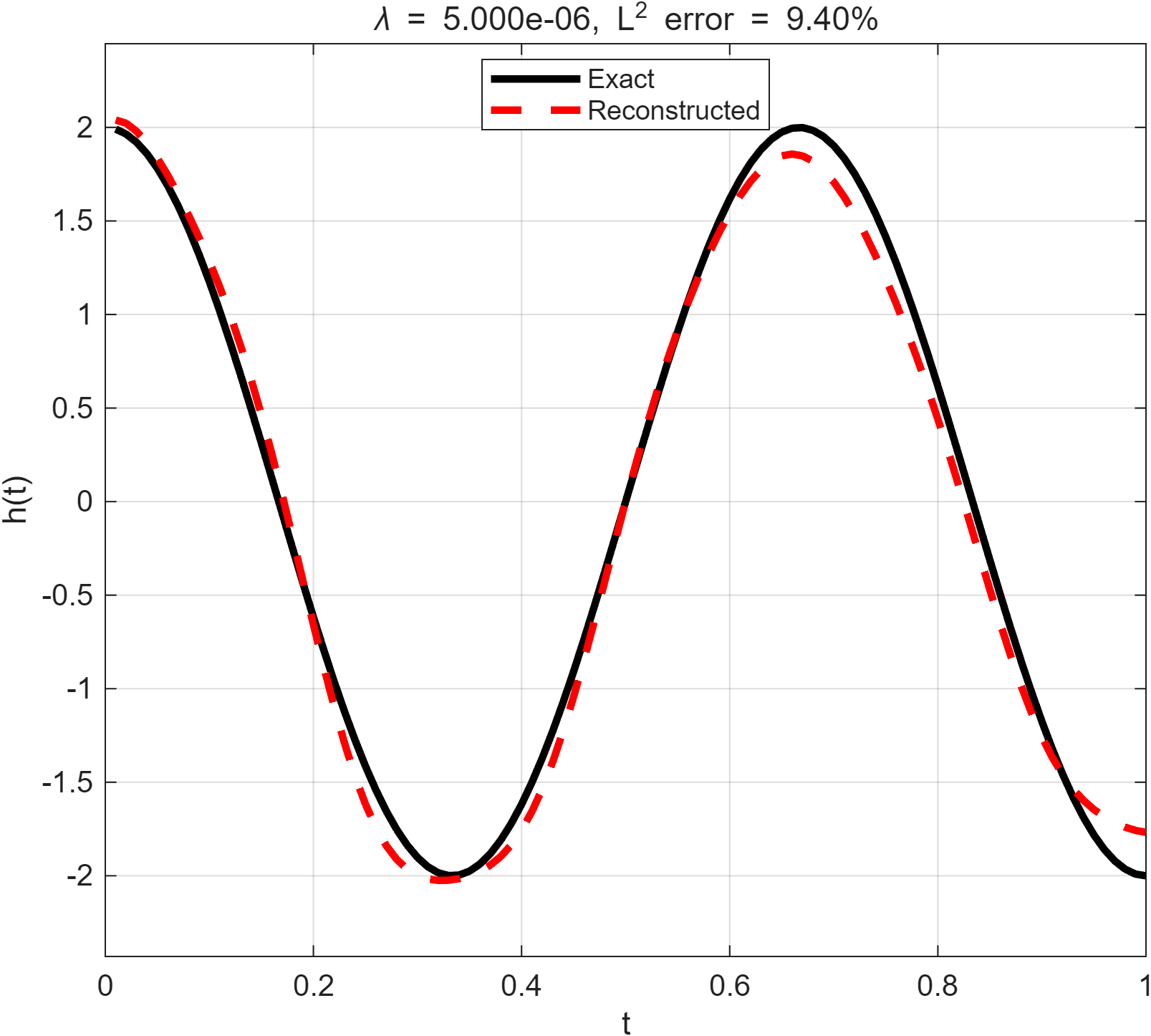}
       \end{minipage}
       \caption{Example \ref{ex2}: Numerical results for the variable fractional order \(\alpha(t) = 0.5 + 0.2t^2\) 
       with \(\delta_{rel} =10\%\): \((\text{left})\,error =6.10\%\), \( (\text{right})\, error=9.40\%\).
 }\label{fig:ex2}
       \end{figure}  
\end{example}

Examples \ref{ex2} and \ref{ex3} are designed to further examine the performance of the proposed
 reconstruction method under different types of model uncertainties. 
 Compared with Example \ref{ex1},  Example \ref{ex2} replaces the linear variable fractional order with a quadratic one, 
 \(\alpha(t)=0.5+0.2t^2\), while keeping all other settings unchanged. 
 As shown in Figure \ref{fig:ex2}, at \(\delta_{rel} = 10\%\), the reconstruction errors are \(6.10\%\) and \(9.40\%\) for the sine and cosine cases, 
 respectively, both slightly smaller than the corresponding errors of \(6.23\%\) and \(9.65\%\) in Example \ref{ex1}. 
 This indicates that the reconstruction accuracy is not highly sensitive to the specific choice of the fractional order function.

\begin{example}\label{ex3}
    Consider the problem \eqref{eq:variable-exponent-subdiffusion} 
    with the variable fractional order \(\alpha(t) = 0.5 + 0.2t\) for $t\in [0,1]$. we test the following cases:
\begin{itemize}
    \item[\bf {(3a):}] \(h(t) = 4t^{\frac{3}{4}}(\frac{1}{4}-t^{\frac{3}{4}} )\) (left in Figure~\ref{fig:ex3});
    \item[\bf {(3b):}] \(h(t) = 1-\left\lvert 2t-1\right\rvert \) (right in Figure~\ref{fig:ex3}).
\end{itemize}

  \begin{figure}[H]
       \centering
       \begin{minipage}[t]{0.48\textwidth}
       \centering
       \includegraphics[width=7.5cm]{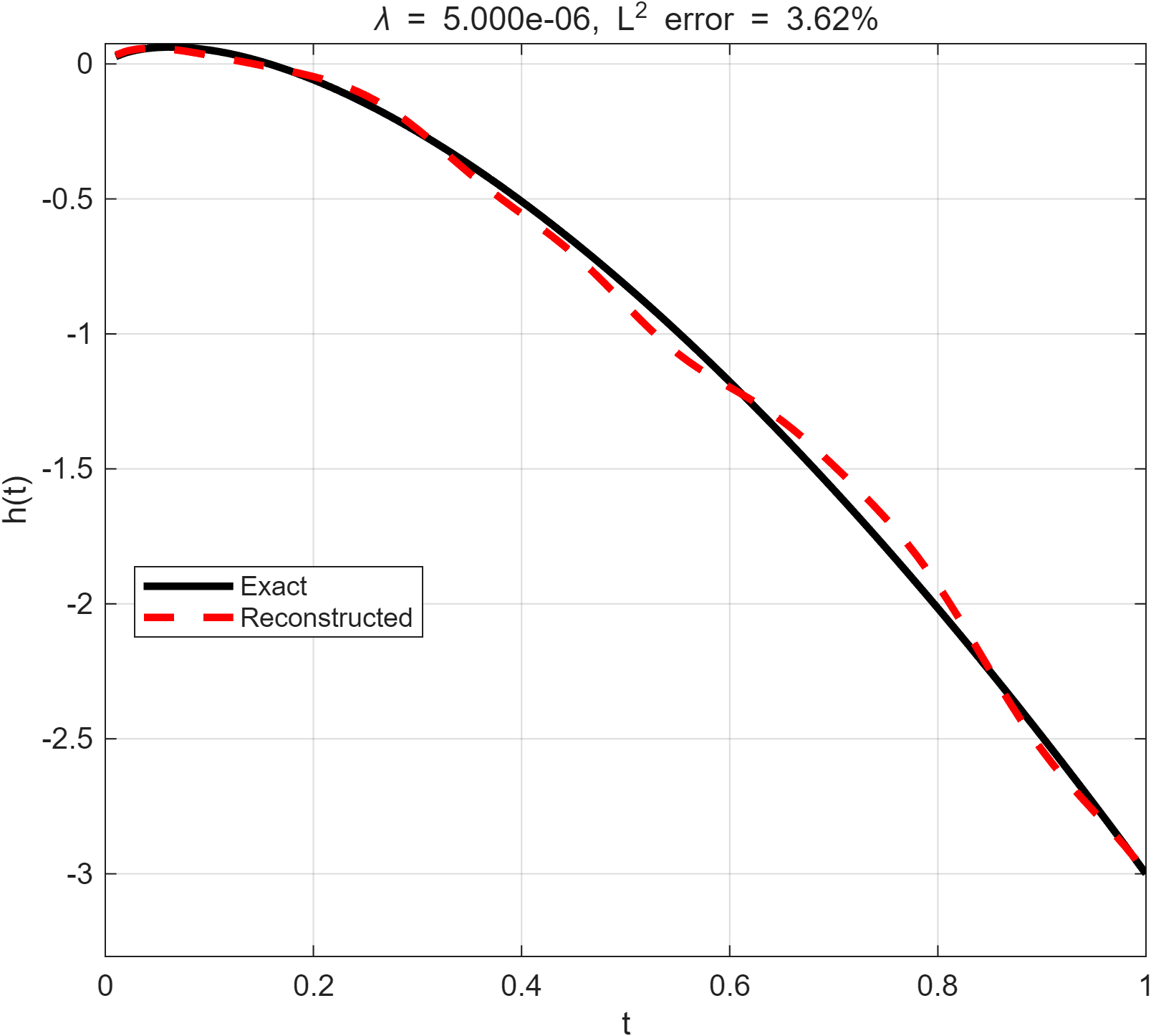}
       \end{minipage}
       \begin{minipage}[t]{0.48\textwidth}
       \centering
       \includegraphics[width=7.5cm]{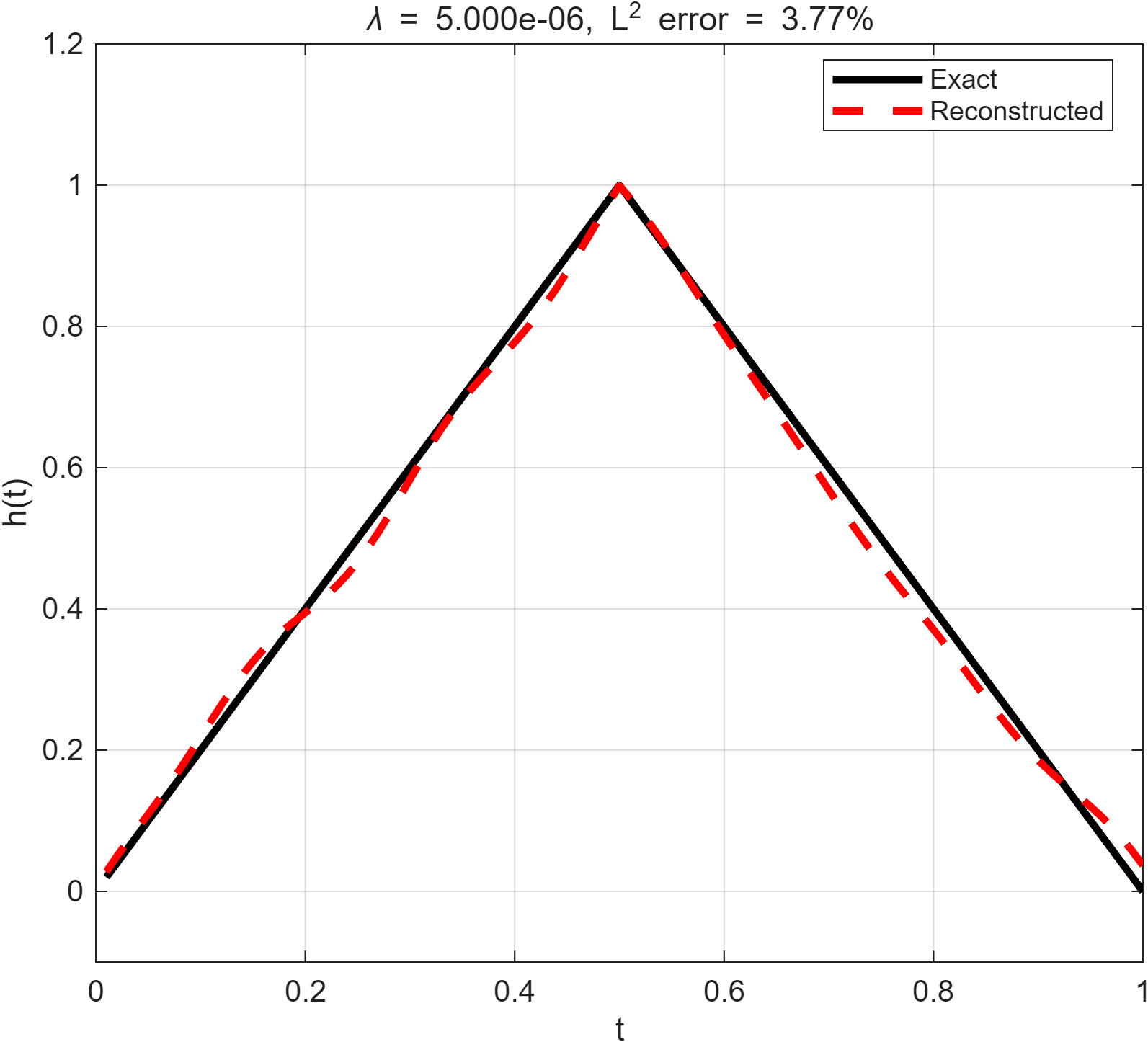}
       \end{minipage}
       \caption{Example \ref{ex3}: Numerical results for the variable fractional order \(\alpha(t) = 0.5 + 0.2t\) 
       with  \(\delta_{rel} =10\%\): \((\text{left})\,error =3.62\%\), \( (\text{right})\, error=3.77\%\).
  }\label{fig:ex3}
       \end{figure}  
\end{example}

In contrast to Examples  \ref{ex1} and \ref{ex2}, Example \ref{ex3} tests two non-smooth source functions: \(h(t)=4t^{3/4}(1/4-t^{3/4})\), 
which exhibits a weak singularity at \(t=0\); and \(h(t)=1-|2t-1|\), which is continuous but not differentiable at \(t=1/2\). 
As shown in Figure \ref{fig:ex3}, despite the lack of smoothness, the proposed method still yields numerically reasonable reconstructions 
at \(\delta_{rel}=10\%\), although some oscillations are observed in the reconstructed curves. Overall, 
these results indicate that the proposed conjugate gradient method with Tikhonov regularization is robust to 
different fractional order functions and can also cope, to some extent, with the lack of temporal regularity in the source term.

\section*{Acknowledgments}

This work was supported by the National Key Research and Development Program of China (Grant No. 2023YFB3002400).

\newpage

\bibliography{vary_bib}

\end{document}